\documentclass[conference]{IEEEtran}
\IEEEoverridecommandlockouts
\usepackage{cite}
\usepackage{amsmath,amssymb,amsfonts}
\usepackage{algorithmic}
\usepackage{graphicx}
\usepackage{hyperref}
\usepackage{textcomp}
\usepackage{xcolor}
\def\BibTeX{{\rm B\kern-.05em{\sc i\kern-.025em b}\kern-.08em
    T\kern-.1667em\lower.7ex\hbox{E}\kern-.125emX}}
\begin{document}

\title{Performance Evaluation of Mixed-Precision Runge-Kutta Methods}

\author{\IEEEauthorblockN{Ben Burnett}
\IEEEauthorblockA{\textit{University of Massachusetts Dartmouth} \\
bburnett@umassd.edu}
\and
\IEEEauthorblockN{Sigal Gottlieb}
\IEEEauthorblockA{\textit{
University of Massachusetts Dartmouth} \\
sgottlieb@umassd.edu
} 
\and
\IEEEauthorblockN{Zachary J. Grant}
\IEEEauthorblockA{\textit{
Oak Ridge National Laboratory} \\
grantzj@ornl.gov}
\and
\IEEEauthorblockN{}
\IEEEauthorblockA{}
\and
\IEEEauthorblockN{Alfa Heryudono}
\IEEEauthorblockA{\textit{University of Massachusetts Dartmouth} \\
aheryudono@umassd.edu}
}

\maketitle

\begin{abstract}
Additive Runge-Kutta methods designed for preserving highly accurate solutions in 
mixed-precision computation were proposed and analyzed in \cite{zack}. 
These specially designed methods use reduced precision for the implicit  computations
and full precision for the explicit computations.
We develop a FORTRAN code to solve a nonlinear system of ordinary differential equations using 
the mixed precision additive Runge-Kutta (MP-ARK) methods on IBM POWER9 and Intel x86\_64 chips.
The convergence, accuracy, runtime, and energy consumption of these methods is explored.
We show that these MP-ARK methods efficiently produce accurate solutions with 
significant reductions in runtime (and by extension energy consumption).
\end{abstract}


\section{Introduction}
 
The use of mixed precision to enhance performance of numerical algorithms
is gaining in popularity, as evidenced by the emergence of mixed precision algorithms
for common linear algebra operations\cite{b1} to aid in the evaluation of large computational 
networks for deep learning.
Mixed precision simulations aim to combine the efficiency of low precision with the accuracy of high precision computations. However, care is needed here to prevent the low precision computations 
from degrading the overall accuracy of the method,
or the high precision computations from adversely impacting the efficiency.

In \cite{perf1}, the authors explored several mixed precision algorithms for solving sparse linear systems and identified algorithms that benefit from a mixed precision implementation and  others that do not. However, while the performance of mixed precision algorithms for linear algebra applications
has been studied \cite{perf1} \cite{perf2} \cite{perf3} \cite{perf4} \cite{perf5} \cite{perf6},
mixed precision algorithms for solving ordinary differential equations (ODES) have only recently been proposed. In \cite{zack} Z. Grant proposed a numerical analysis framework for development of mixed-precision Runge--Kutta methods by using an additive Runge--Kutta method approach and framing the use of multiple  precisions as a perturbation. Using this rigorous approach, \cite{zack} developed novel methods that reduce the cost of the  computationally expensive implicit stages in the Runge--Kutta methods  by employing a low precision computation. Meanwhile, the structure of the MP-ARK is designed to suppress the low precision errors either by introducing inexpensive explicit high-precision correction terms,
or by designing novel methods that internally suppress the low precision perturbations.

The work in \cite{zack} focused on the development of a framework for analysis of MP-ARK methods, and
the development and convergence verification of such MP-ARK methods, using a simplified "chopping"
routine to simulate low precision in MATLAB.
In this work, we revisit the novel methods in \cite{zack} and evaluate their performance on multiple CPU architectures in single, double, and quadruple precision. 
We first confirm that these methods
converge as expected. Next, we demonstrate that for a given level of accuracy, the mixed precision
codes provide reductions in runtime (and by extension energy consumption) ranging from 2 
times speed up to upwards of 15 times speed up.

\section{MP-ARK methods} \label{methods}

For this study we looked at three of the methods outlined in \cite{zack}: the mixed precision implicit midpoint method
with and without correction steps, the mixed-precision SDIRK method with and without correction steps, 
and a novel MP-ARK method from \cite{zack}
designed to suppress the low precision errors without additional correction steps. 
In this section we describe  these methods.

\subsection{Mixed Precision Methods Based on the Implicit Midpoint Rule}
The implicit midpoint method was the main motivating method of the original paper. 
It was used as the base example for  mixed precision methods 
and a demonstration of the perturbed Runge-Kutta analysis to them. 
The base method is given by
\begin{eqnarray*}\label{midpoint}
y^{(1)} &=& u^n + \frac{\Delta t}{2} F\big(y^{(1)}\big) \\
 u^{n+1} &=& u^n + \Delta t F\big(y^{(1)}\big)
\end{eqnarray*}
The  mixed precision formulation of this method is then defined as
\begin{subequations} \label{MPmidpoint1}
\begin{eqnarray}
y^{(1)} &=& u^n + \frac{\Delta t}{2} F^{\epsilon}\big(y^{(1)}\big)  \\ 
 u^{n+1} &=& u^n + \Delta t F\big(y^{(1)}\big)
\end{eqnarray}
\end{subequations}
where $\epsilon$ represents the precision level for the low precision evaluation.

The low precision evaluation of the implicit solver in the first stage of Equation \eqref{MPmidpoint1} 
requires a subroutine that maps the variables low precision, performs a Newton iteration
with tolerances that are consistent with the level of the low precision, and 
returns the low precision value of the variable $y^{(1)}$. This value is then cast to
a high precision variable $y^{(1)}$ to be used in the explicit second step of 
 Equation \eqref{MPmidpoint1}. 
 
Adding explicit high precision correction steps to the mixed precision implicit midpoint rule 
\eqref{MPmidpoint1} we define the {\em corrected MP-IMR} with $p-1$ correction terms
\begin{subequations} \label{MPmidpoint2}
\begin{eqnarray}
 y_{[0]}^{(1)} & = & u^n + \frac{\Delta t}{2} F^{\epsilon}\big(y_{[0]}^{(1)}\big) \\
 y_{[k]}^{(1)} & = &  u^{n} + \frac{\Delta t}{2} F\big(y_{[k-1]}^{(1)}\big) \; \text{for}\; k=1,...,p-1  \\
 u^{n+1} & = & u^n + \Delta t F\big(y_{[p-1]}^{(1)}\big)
\end{eqnarray}
\end{subequations}
The additional correction stages  are performed at high precision  but they are cheap explicit computations, 
while the implicit solves, though performed in low precision, 
are expected to be the computationally dominant.
The  rationale behind the correction stage approach is that the gain in accuracy will be well-worth the extra two cheap explicit stage. Of course, the purpose of our numerical experiments will be to verify this assumption.

\subsection{Mixed Precision SDIRK method with and without corrections}

The singly diagonally implicit Runge-Kutta (SDIRK) method is a higher order than the implicit midpoint method and involves two implicit stages. The mixed precision SDIRK method with $p-1$ correction terms is  (where $k=1:p-1$)
\begin{eqnarray}\label{MP-SDIRKc}
 y_{[0]}^{(1)} & = &  u^n + \gamma \Delta t F^{\epsilon}\big(y_{[0]}^{(1)}\big) \nonumber \\
 y_{[k]}^{(1)} & = & u^n + \gamma \Delta t F\big(y_{[k-1]}^{(1)}\big)   \\
 y_{[0]}^{(2)} & = & y^n + (1 - 2 \gamma)
 \Delta t F \big(y_{[p-1]}^{(1)}\big) 
+ \gamma \Delta t F^{\epsilon}\big(y_{[0]}^{(2)}\big)  \nonumber \\
 y_{[k]}^{(2)} & = & y^n + (1 - 2 \gamma)\Delta t F\big(y_{[p-1]}^{(1)}\big) 
 + \gamma \Delta t F\big(y_{[k-1]}^{(2)}\big) 
 \nonumber \\
 y^{n+1} & = & u^n + \frac{\Delta t}{2} F\big(y_{[p-1]}^{(1)}\big) + 
 \frac{\Delta t}{2} F\big(y_{[p-1]}^{(2)}\big) .
 \nonumber
\end{eqnarray}

\subsection{Mixed precision 4s3p MP-ARK method}

The novel mixed precision additive Runge-Kutta method (MP-ARK) developed by Z. Grant \cite{zack} is defined by

\begin{subequations}
\begin{eqnarray}\label{MP-4s3pA}
y^{(1)} & = & u^n + \Delta t A^{\epsilon}_{1,1} F^{\epsilon} \big(y^{(1)}\big) \\
y^{(2)} & = & u^n + \Delta t A_{2,1} F \big(y^{(1)}\big) \\
y^{(3)} & = & u^n + \Delta t \big[ A_{3,1} F \big(y^{(1)}\big) + A_{3,2} F \big(y^{(2)}\big) \big] \\
 &+& \Delta t \big[ A^{\epsilon}_{3,1} F^{\epsilon} \big(y^{(1)}\big) + A^{\epsilon}_{3,3} F^{\epsilon} \big(y^{(3)}\big) \big] \nonumber \\
y^{(4)} & = & u^n + \Delta t \big[ A_{4,1} F \big(y_{(1)}\big) + A_{4,2} F \big(y^{(2)}\big) \\ 
 &+& A_{4,3} F \big(y^{(3)}\big) \big] \nonumber \\
u^{n+1} & = & u^{n} + \frac{1}{2}\Delta t \big[ \big(y^{(2)}\big) + \big(y^{(4)}\big) \big]
\end{eqnarray}
\end{subequations}
where $A$ and $A^{\epsilon}$ are 
{\small
\begin{eqnarray*}\label{butcher_coeffs}
A_{2, 1}  =  0.211324865405187, &
A_{3, 1}  =  0.709495523817170,  \\
A_{3, 2}  =  -0.865314250619423, &
A_{4, 1}  =  0.705123240545107,   \\
A_{4, 2}  = 0.943370088535775, &
A_{4, 3} =  -0.859818194486069,   \\
A^{\epsilon}_{1, 1}  =  0.788675134594813, &
A^{\epsilon}_{3, 1}  =  0.051944240459852, \\
A^{\epsilon}_{3, 3} =  0.788675134594813. &
\nonumber  
\end{eqnarray*}
}
\vspace{-.15in}
\section{Implementation Details}

The mixed-precision Runge-Kutta methods that we explored were implemented in Fortran using the 2008 standard. 
The \verb|iso_fortran_env| intrinsic module was also used to import the \verb|Real32|, \verb|Real64|, and \verb|Real128| derived types.
The implicit solvers were implemented using the Newton-Raphson method for multivariate systems. 

As in \cite{zack}, we solve the Van der Pol system 
\begin{subequations} \label{vdp}
\begin{eqnarray}
y_1' & = & y_2 \\
y_2' & = &  y_2 (1 - y_1^2) - y_1
\end{eqnarray}
\end{subequations}
 with  initial conditions $y_1(0) = 2.0$, $y_2(0) = 0.0$, for time $t = [0,1]$. A reference solution using an explicit fourth order Runge--Kutta at a small time step was used for calculating the errors. The reference solution was computed entirely in quadruple precision and all methods results regardless of computational precision, were cast to quadruple precision in order to compute the error between a method and this reference solution.

The tolerance of the implicit solvers were set using the Fortran intrinsic machine epsilon $\epsilon$, where  the tolerance value  $(1 + 0.001)\epsilon$ was chosen. 
The maximum number of iterations for the solver was set to $20$ regardless of the precision, however this fail-safe maximum was not attained. 
These setting do have significant impact on the overall runtime of the method and were determined to be fair for the consideration of mixed precision methods, based on our tests showing that the implicit solvers on average only used $4$ iterations to converge.

Higher precision computations  (specifically quadruple precision but also arbitrary precision) 
are generally evaluated in software on a majority of platforms. 
Notable exceptions for this are the new IBM POWER9 
that have hardware level support for quadruple precision. In this work, we compare the runtime performance
on both software based high precision computations and hardware based high precision computations.
The experiments were performed on a system with an Intel Xeon Gold 6126 CPU running CentOS 7 and an IBM POWER9 running CentOS 8. 
We will refer to these as x86 and POWER9, respectively.
The experiment performed looks at three levels of precision: single, double, and quadruple. 
The Intel system is limited to hardware evaluation for single and double precision and requires 
quadruple precision to be evaluated using software. The POWER9 system is capable of hardware 
evaluation of all precisions. This was done to show the usefulness of these methods under 
both conditions and that the benefits are platform independent. 

The codes were compiled on both architectures using gfortran version 11.1.0 using no optimization flags. 
On the POWER9 system, the GCC compilers were compiled using the IBM AT14.0 libraries 
in order to enable hardware quadruple precision support while the GCC compilers on the Intel system
were compiled using the default compiler settings. Timing on the method subroutine was done using the 
intrinsic \verb|cpu_time| Fortran subroutine. 

The code for this is available on GitLab at \url{https://gitlab.com/bburnett6/mixed-precision-rk}. 
This repository contains all the necessary files to reproduce the study performed here. 
The only dependencies are a working version of the GCC 11.1.0 compilers and an environment 
with the Numpy and MatPlotLib python libraries for the experiment driver and analysis 
files that were implemented in python. 

\section{Results}
 In this section we describe the errors, runtime, and energy consumed when using the different 
 approaches described in Section \ref{methods} to solve the non-stiff Van der Pol system \eqref{vdp}.

\subsection{Mixed precision implicit midpoint method}
\subsubsection{Accuracy study}
\begin{figure}[htbp]
\centerline{\includegraphics[width=0.5\textwidth]{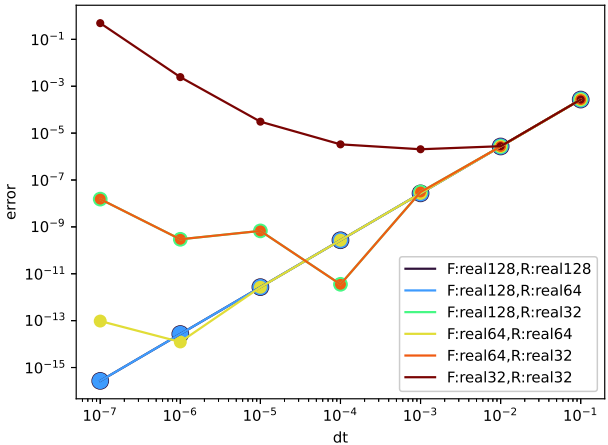}}
\caption{Errors from mixed precision implicit midpoint rule}
\label{impmid_err}
\end{figure}

We begin with the mixed precision implicit midpoint rule \eqref{MPmidpoint1}. 
We reproduced the convergence behavior reported in \cite{zack}, as we show in Figure \ref{impmid_err}. 
We observe that the single precision
computation is not convergent for many values of  $\Delta t$: the error builds up due to roundoff for smaller $\Delta t$.
The mixed precision computations that use single precision (double-single or quad-single)
perform better, reducing the errors to below $10^{-7}$; however, they still show evidence of degredation of accuracy due to the single precision computations. The mixed precision quad-double computation performs as well as the quad-quad computation for all values of $\Delta t$ tested. However, the double precision computation performs just as well for all but the smallest $\Delta t$ tested.

\begin{figure}[htbp]
\centerline{\includegraphics[width=0.5\textwidth]{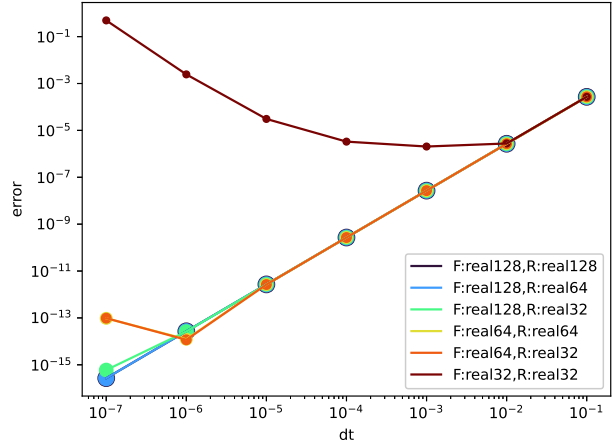}}
\caption{Errors from mixed precision implicit midpoint rule with one correction step}
\label{impmid_1p_err}
\end{figure}

\begin{figure}[htb]
\includegraphics[width=0.5\textwidth]{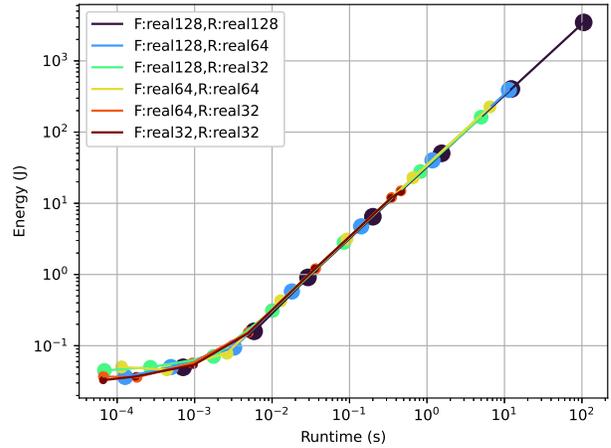}
  \caption{Total energy used by the mixed precision implicit midpoint rule with 
  one correction scales linearly with the runtime.}
  \label{energy-runtime}
  \end{figure}
  

The corrected mixed precision implicit midpoint rule \eqref{MPmidpoint2} 
with only one correction fixes the problems observed in the mixed precision 
computation of \eqref{MPmidpoint1}. As we observe in Figure \ref{impmid_1p_err}, 
all the mixed precision computations have the same errors as their high-precision 
counterparts. An additional correction (not 
shown) enhances the errors in the mixed precision simulation a bit further.

\subsubsection{Runtime \& energy consumption}
\begin{figure*}[htbp]
\includegraphics[width=0.495\textwidth]{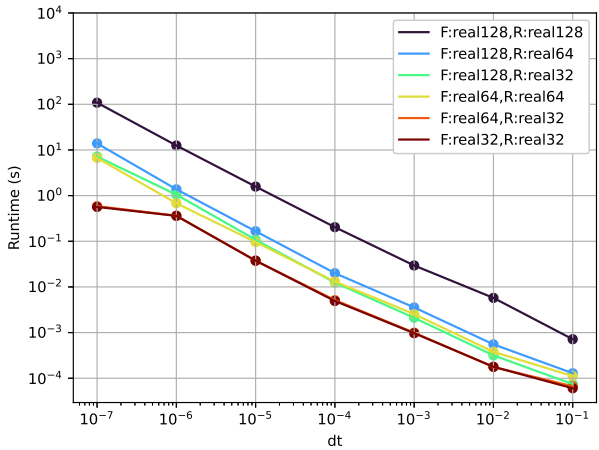}
\includegraphics[width=0.495\textwidth]{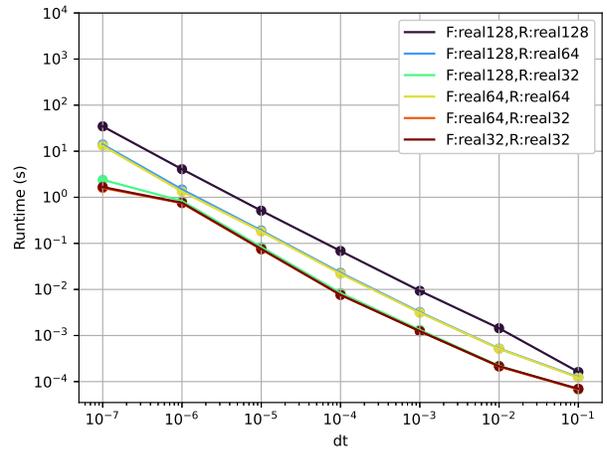}
\caption{Mixed precision implicit midpoint method with two corrections: Time to Completion. Left: x86,
Right: POWER9.}
\label{impmid_timeP2}
\end{figure*}

The energy consumption results for this simulation follow very closely the runtime results.
  Figure \ref{energy-runtime} shows that for long enough runtimes, the energy consumed scales 
  linearly with the runtime for all precisions.
In Figure  \ref{impmid_timeP2} 
we present the runtime of the mixed precision implicit midpoint rule with two corrections,
on the x86 and POWER9 chips, respectively. The runtime is not impacted by the correction terms. 
\begin{figure*}[htb]
  \includegraphics[width=0.495\textwidth]{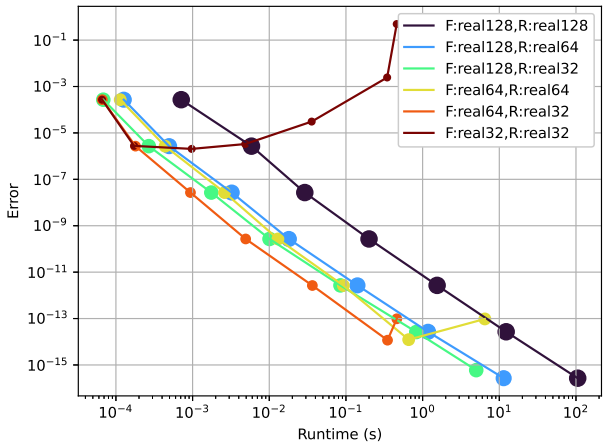}
   \includegraphics[width=0.495\textwidth]{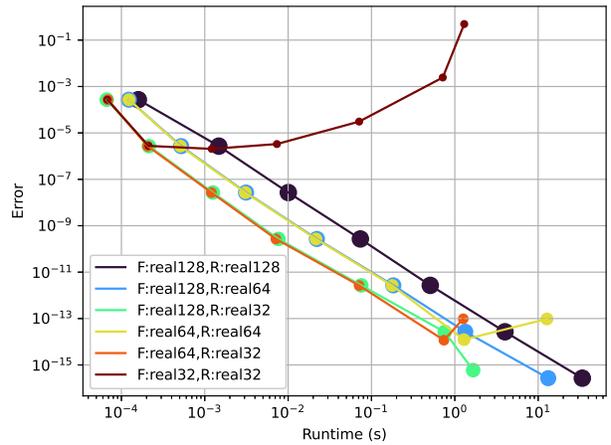}
\caption{Mixed precision implicit midpoint method with one correction step: error vs runtime. Left: x86; right: POWER9}
\label{impmid_1p_time}
\end{figure*}

Remarkable in this figure is the significant reduction in runtime 
from the quad precision to the mixed quad-double precision for the x86 chip, 
and from the quad precision to the mixed quad-single precision for the POWER9 chip. 
The cost of quad precision on the POWER9 chip is significantly less than on the x86 chips,
but the runtime savings is not  impacted, as we see in the tables below.
The differences in the runtime behaviors is explained by the fact that POWER9 chip
has hardware-level support for quad precision.

\subsubsection{Efficiency}
To better illustrate the benefit of the mixed precision computation we plot 
the method error against the runtime for the mixed precision implicit midpoint method with
one correction, in Figure \ref{impmid_1p_time}.

\bigskip
\begin{tabular}{ |c|c|c|c|c|c|  }
 \hline
 \multicolumn{6}{|c|}{Mixed precision IMR method with two corrections}\\
  \multicolumn{6}{|c|}{Runtime (s) for x86} \\ \hline
 error & 64/64 & 64/32 & 128/128 & 128/64 
 & 128/32 \\
$\approx 10^{-9}$ & 0.013 & 0.005 & 0.203 & 0.020 & 0.012 \\
$\approx 10^{-15}$ & N/A & N/A & 108.4 & 13.91 & 7.165 \\ \hline
 \end{tabular}
 \bigskip
 
On the x86 (Figure \ref{impmid_1p_time}, left), we see that the mixed precision double-single code is the {\em most} efficient. For an error level of 
$\approx 10^{-15}$, the mixed precision quad-single code has a fifteen-fold reduction in runtime.
The mixed precision double-single code gives the same level of error at $\approx 10^{-9}$
as the double precision code with a runtime reduction of 2.6-fold.

On the POWER9 (Figure \ref{impmid_1p_time}, right; also table below) the runtime is generally better, but the
runtime   savings are similar. For the  error level of  $\approx 10^{-9}$, we see 3-fold savings from the mixed precision
double-single over the the double precision. For the  error level of  $\approx 10^{-15}$, 
we see more than 14-fold savings from the mixed precision quad-single over the the quad precision. 

\smallskip
\begin{tabular}{ |c|c|c|c|c|c|  }
 \hline
 \multicolumn{6}{|c|}{Mixed precision IMR method with two corrections}\\
  \multicolumn{6}{|c|}{Runtime (s) for POWER9} \\ \hline
 error & 64/64 & 64/32 & 128/128 & 128/64 
 & 128/32 \\
$\approx 10^{-9}$ & 0.021 & 0.007 & 0.0683 & 0.023 & 0.008 \\
$\approx 10^{-15}$ & N/A & N/A & 34.64 & 14.00 & 2.395 \\ \hline
 \end{tabular}

\begin{figure*}[ht!]
  \includegraphics[width=0.5\textwidth]{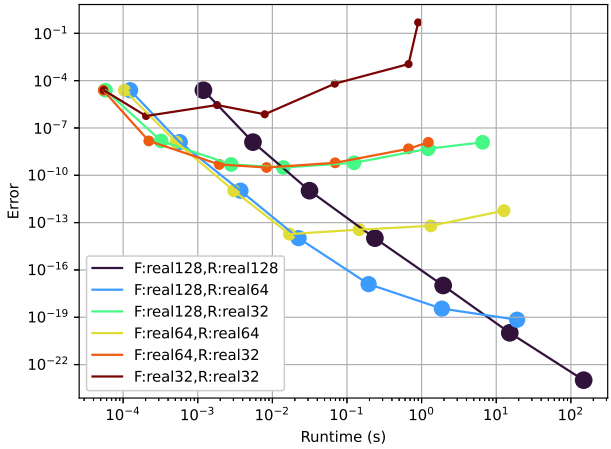}
  \includegraphics[width=0.5\textwidth]{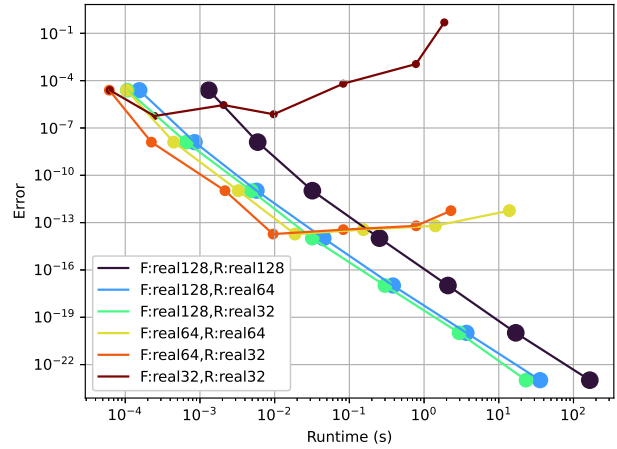}\\
  \caption{Mixed precision SDIRK method errors vs runtime on x86. 
  Top left: no corrections;
  top right: three corrections.}
  \label{sdirk-eff-x86}
\end{figure*}

\subsection{Mixed precision SDIRK Method}
\begin{figure}[h!]
  \begin{center}
  \includegraphics[width=0.225\textwidth]{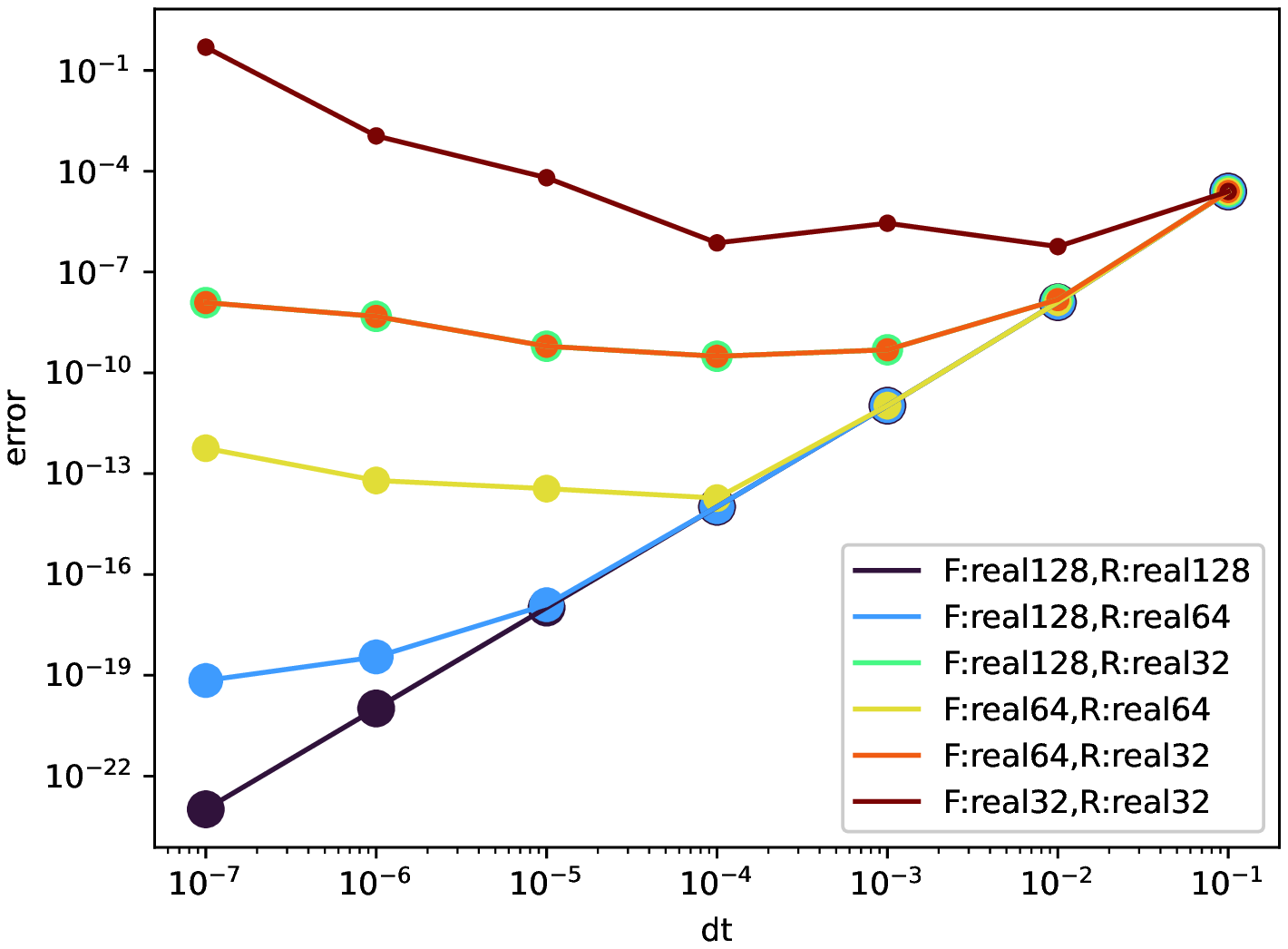}
  \includegraphics[width=0.225\textwidth]{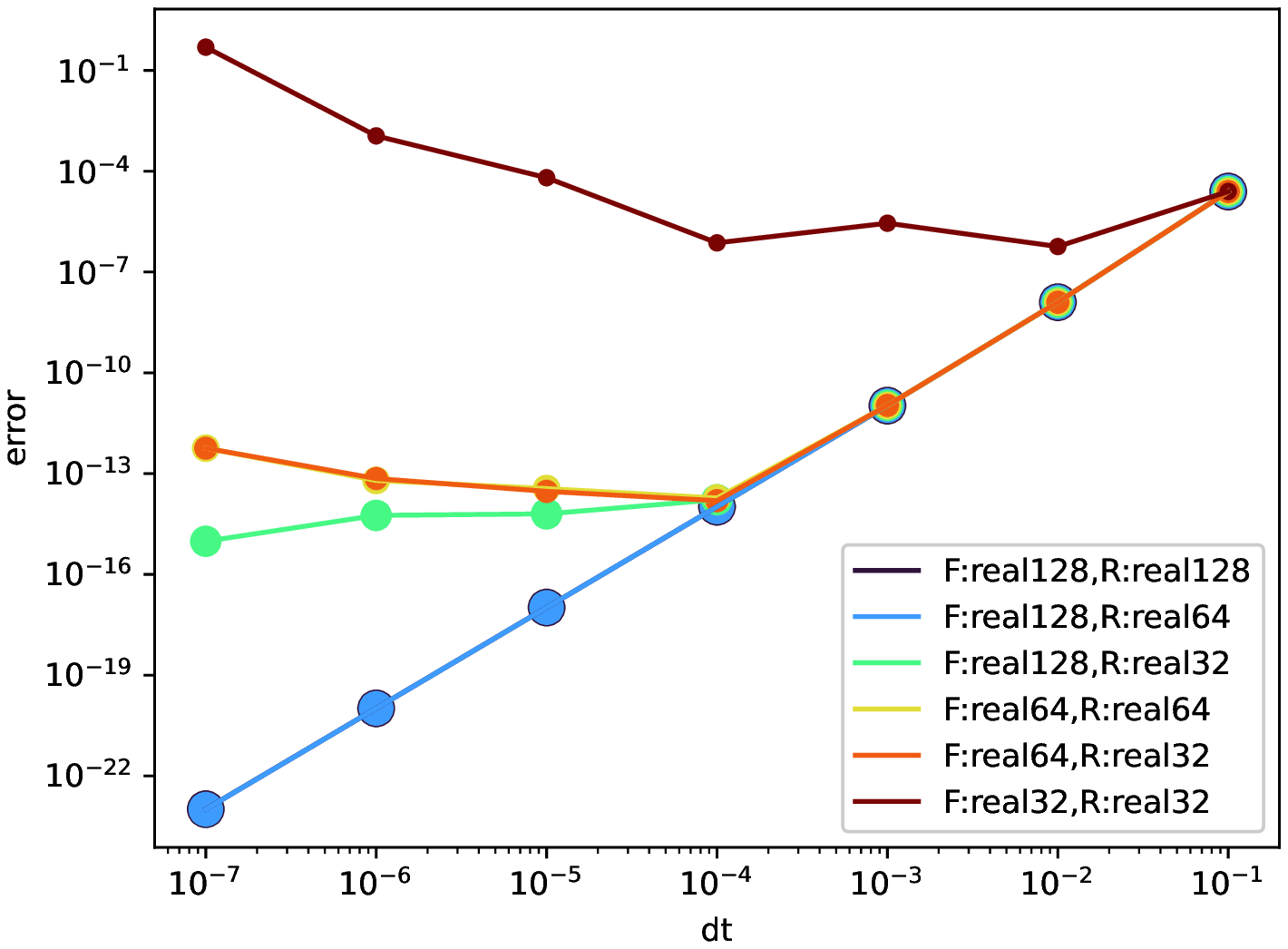} \\
  \includegraphics[width=0.225\textwidth]{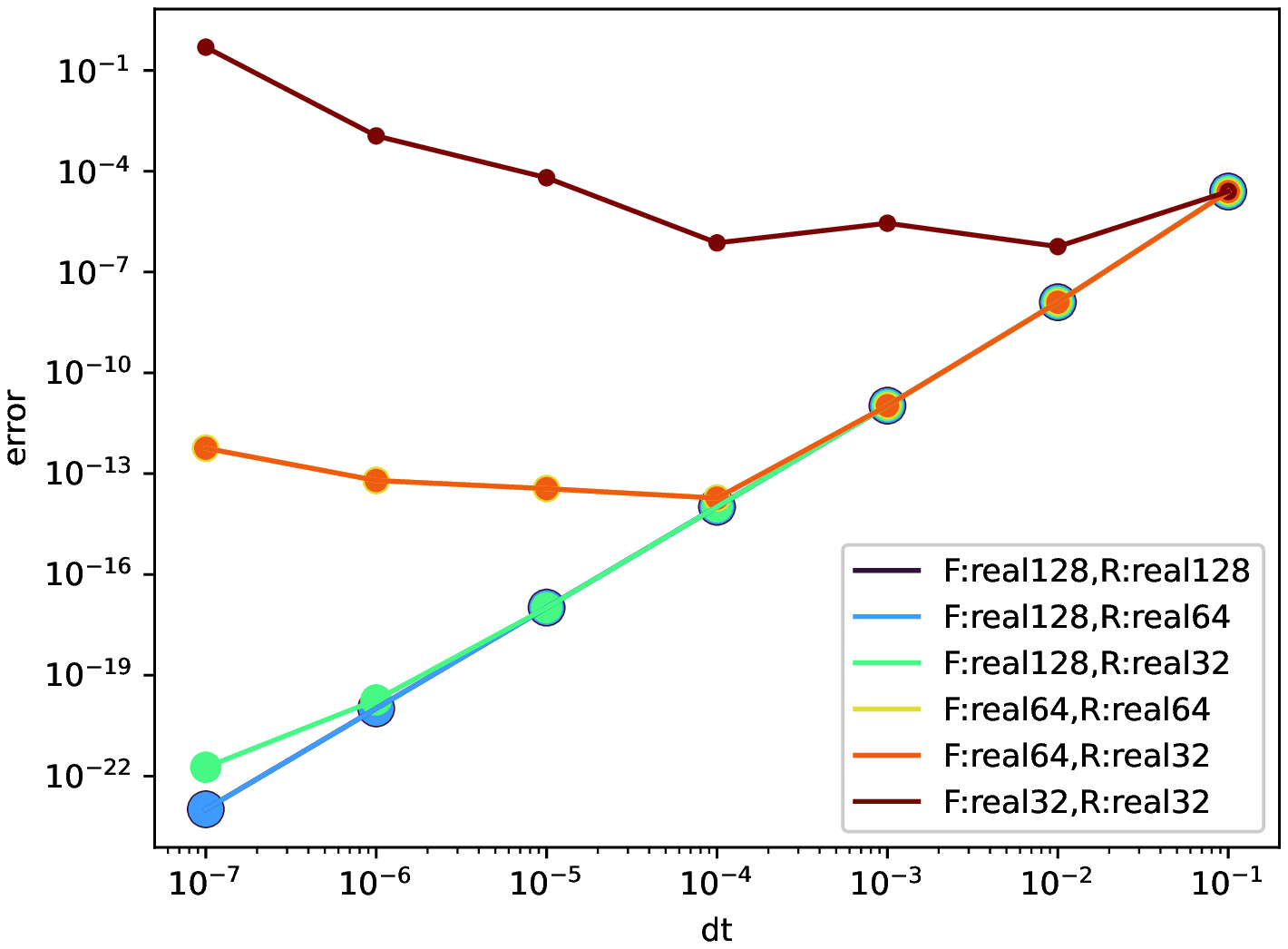}
  \includegraphics[width=0.225\textwidth]{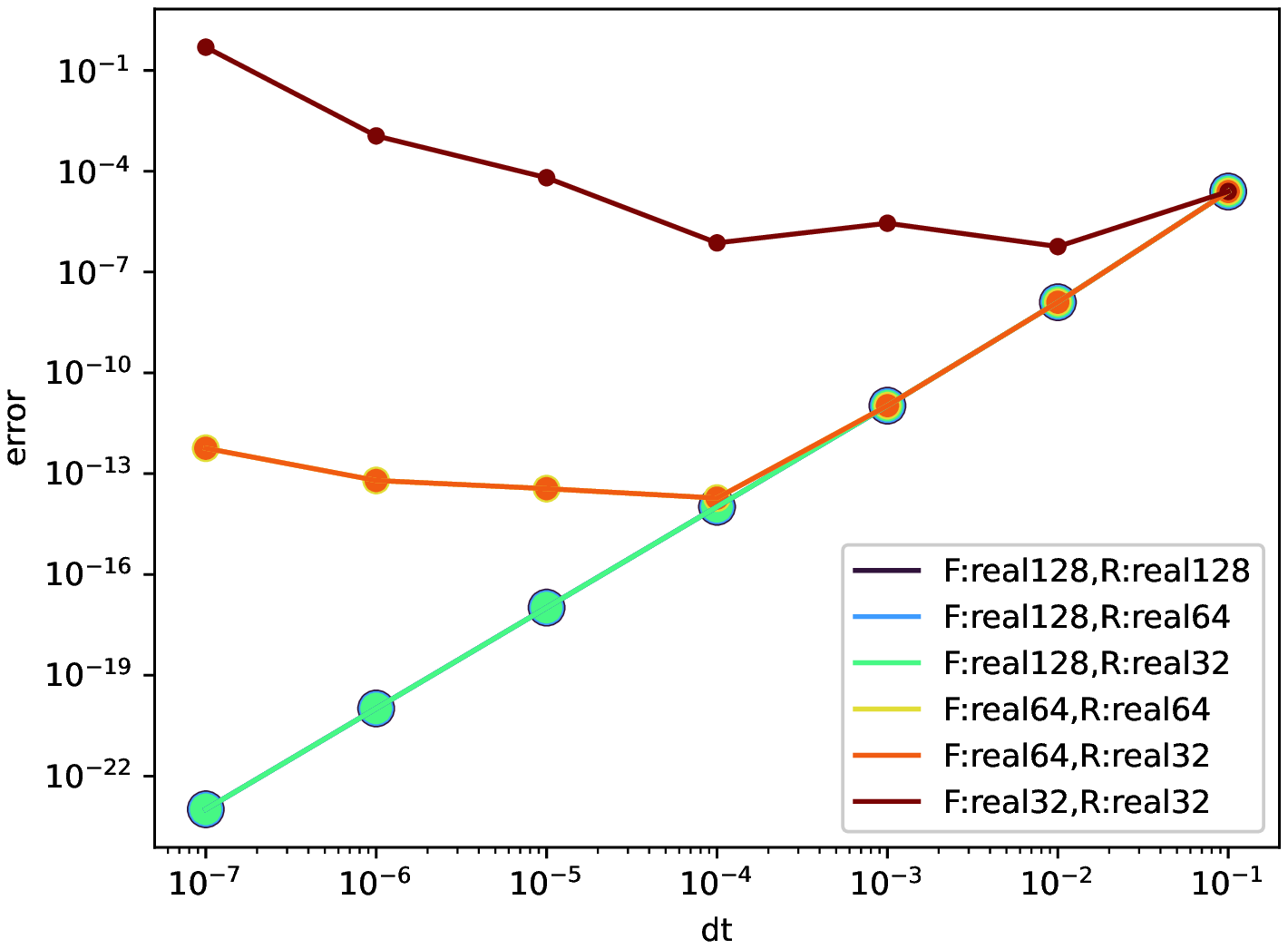}\\
  \end{center}
  \caption{Mixed precision SDIRK method errors. 
Top left: no corrections;
Top right: one correction;
Bottom left: two corrections;
Bottom right: three corrections.
  }
  \label{sdirk_err}
\end{figure}

The results for the mixed precision singly diagonally implicit Runge-Kutta (MP-SDIRK) in 
Section \ref{methods} are similar to those of the implicit midpoint rule.
Figure \ref{sdirk_err} shows the convergence of this method with no corrections, 
one correction, two corrections, and three corrections. The results show clearly 
how with each successive correction the errors from the  mixed precision simulation are reduced, 
until the mixed precision results become as accurate as their full-precision counterparts. 
Furthermore, we confirm that additional corrections 
do not result in a detrimental increase in runtime. 

In Figure \ref{sdirk-eff-x86} we show the errors per runtime for the mixed precision SDIRK method with no corrections (left) and with three corrections (right), on the x86 chip. We observe that the method with three corrections is overall significantly more efficient than the method with no corrections. 
For errors  above $\approx 10^{-13}$, the mixed precision double-single method performs best. 
For errors below that level, the mixed precision 
 quad-single method is most efficient, as shown in the table below.
  The runtime savings on the POWER9 chip are similar.
 
 \smallskip
 
 \begin{tabular}{ |c|c|c|c|c|c|  }
 \hline
 \multicolumn{6}{|c|}{Mixed precision SDIRK with three corrections}\\
  \multicolumn{6}{|c|}{Runtime (s) for x86} \\ \hline
 error & 64/64 & 64/32 & 128/128 & 128/64 
 & 128/32 \\
$\approx 10^{-13}$ & 0.018 & 0.009 & 0.253 & 0.045 & 0.031 \\
$\approx 10^{-22}$ & N/A & N/A & 165.3 & 35.72 & 23.15 \\ \hline
 \multicolumn{6}{|c|}{Runtime (s) for POWER9} \\ \hline
 error & 64/64 & 64/32 & 128/128 & 128/64 
 & 128/32 \\
$\approx 10^{-13}$ & 0.031 & 0.016 & 0.092 & 0.034 & 0.019\\
$\approx 10^{-22}$ & N/A & N/A & 54.33 & 30.06 & 8.151\\ \hline
 \end{tabular}

\smallskip

\subsection{Novel Method 4s3pA}
\begin{figure}[htbp]
\centerline{\includegraphics[width=0.5\textwidth]{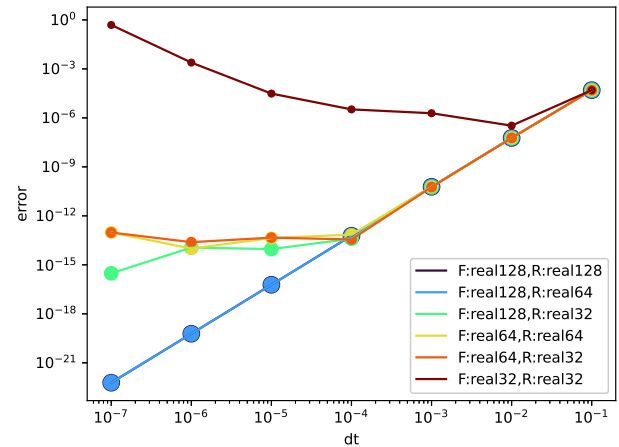}}
\caption{Novel mixed-precision method 4s3pA Method Error}
\label{novelA_err}
\end{figure}

Finally, we present the results from the novel four stage third order method presented in \cite{zack} as
4s3pA. The errors per $\Delta t$ are presented in Figure \ref{novelA_err}. 
We observe that the mixed precision double-single has essentially 
the same errors as the fully double precision computation, 
at almost 3-fold savings in runtime.
The quad-double has the same errors as the fully quad precision 
computation at runtime savings factors of 2.6 to 5.75  depending 
on the platform.

\begin{figure}[htbp]
\centerline{\includegraphics[width=0.5\textwidth]{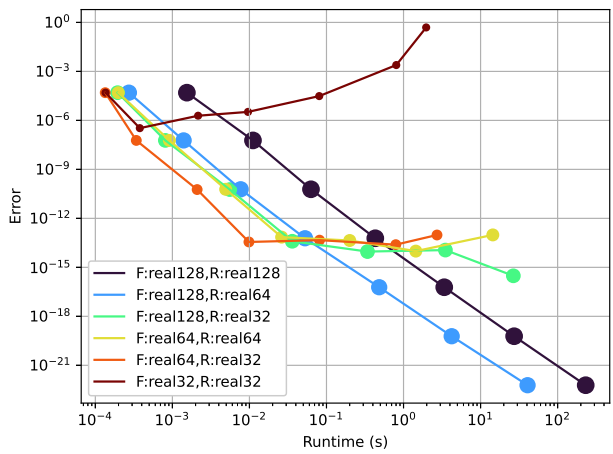}}
\caption{Novel mixed-precision  method 4s3pA Error vs runtime}
\label{novelA_time}
\end{figure}

 \smallskip

\begin{tabular}{ |c|c|c|c|c|c|  }
 \hline
 \multicolumn{6}{|c|}{Mixed precision novel 4s3pA method}\\ \hline
  \multicolumn{6}{|c|}{Runtime (s) for x86} \\ \hline
 error & 64/64 & 64/32 & 128/128 & 128/64 
 & 128/32 \\
$\approx 10^{-13}$ & 0.026 & 0.009 & 0.4300 & 0.052 & 0.035\\
$\approx 10^{-22}$ & N/A & N/A & 231.4 & 40.71 & N/A\\ \hline
\multicolumn{6}{|c|}{Runtime (s) for POWER9} \\ \hline
 error & 64/64 & 64/32 & 128/128 & 128/64 
 & 128/32 \\
$\approx 10^{-13}$ & 0.047 & 0.017 & 0.140 & 0.049 & 0.018\\
$\approx 10^{-22}$ & N/A & N/A & 74.68 & 29.42 & N/A\\ \hline
 \end{tabular}
 
 \smallskip

\section{Conclusions}

In this paper we have presented a numerical study of the new mixed precision Runge--Kutta methods on
two different platforms.
The results presented confirm that the 
mixed precision Runge-Kutta methods on a system
of nonlinear ODEs provide significant runtime (and energy) savings per level of error. This is of particular interest when a high level of accuracy is required, so that expensive high precision must be used. The nonlinear ODE system explored here is a small and relatively simple problem; we expect that these savings would be extended (and indeed improved) on more complicated problems.

\section*{Acknowledgment}
Z. Grant acknowledges the support of the Department of Computational and Applied Mathematics, 
Oak Ridge National Laboratory, Oak Ridge TN 37830. This material is based upon work supported 
by the U.S. Department of Energy, Office of Science, Office of Advanced Scientific Computing Research, 
as part of their Applied Mathematics Research Program. The work was performed at the Oak Ridge National 
Laboratory, which is managed by UT-Battelle, LLC under Contract No. De- AC05-00OR22725. 
The United States Government retains and the publisher, by accepting the article for publication, 
acknowledges that the United States Government retains a non-exclusive, paid-up, irrevocable, 
world-wide license to publish or reproduce the published form of this manuscript, or allow others to do so, 
for the United States Government purposes. The Department of Energy will provide public access to 
these results of federally sponsored research in accordance with the DOE Public Access Plan 
\url{http://energy.gov/downloads/doe-public-access-plan}.

Gottlieb's work was partially supported by AFOSR grant FA9550-18-1-0383.
Burnett, Gottlieb, and Heryudono's work were partially supported by ONR UMass Dartmouth Marine and UnderSea Technology (MUST) grant N00014-20-1-2849 under the project S31320000049160.

All computations for this project were performed on 
CARNIE, the UMass-Dartmouth HPC cluster purchased by
ONR DURIP grant N00014-18-1-2255.

We are grateful to the  IBM Linux team that extended a warm and helping hand to figure out how our POWER9 system was running, and what will be needed to get hardware quadruple precision support.

\clearpage

\end{document}